\input amstex
\documentstyle{amsppt}
\magnification=\magstep1
\vsize 21 true cm
\hsize 14.5 true cm
\abovedisplayskip=1in
\belowdisplayskip=1in
\hoffset=0.5in
\voffset=0.5in

\catcode `\@=11
\let \logo@=\relax
\catcode `\@=\active
\def\dim{\operatorname{dim}}%
\def\depth{\operatorname{depth}}%
\def\rank{\operatorname{rank}}%
\def\Hom{\operatorname{Hom}}%
\def\SHom{\operatorname{\Cal H\text{\it om}}}%
\def\Sing{\operatorname{Sing}}%
\def\reg{\operatorname{reg}}%
\def\codim{\operatorname{codim}}%
\def\Pic{\operatorname{Pic}}%
\def\height{\operatorname{height}}%
\def\Ker{\operatorname{Ker}}%
\def\Proj{\operatorname{Proj}}%

\topmatter
\title\nofrills
{Castelnuovo-Mumford regularity bound for smooth threefolds
in $\Bbb P^5$ and extremal examples}
\endtitle
\leftheadtext\nofrills{\it S. Kwak, Castelnuovo-Mumford
regularity bound}
\rightheadtext\nofrills{\it S. Kwak, Castelnuovo-Mumford
regularity bound}
\author
{\rm By} {\it Si-Jong Kwak} {\rm at KIAS}
\endauthor
\address
School of  Mathematics,  Korea Institute for Advanced Study, Seoul,
Korea 
\endaddress
\email
kwak\@math.snu.ac.kr, sjkwak\@kias.kaist.ac.kr
\endemail

\thanks
This work is partially supported by KIAS and GARC(Global Analysis
Research Center), Seoul, Korea
\endthanks

\abstract
Let  $X$  be a nondegenerate integral subscheme of  dimension
$n$ and  degree  $d$  in  ${\Bbb P}^N$ defined over 
the complex number field ${\Bbb C}$.  
$X$  is said to be $k$-regular  if  
$H^i ({\Bbb  P}^N , {\Cal  I}_X(k-i))  =  0  \quad  \text{for  all}
\,\,   i  \ge  1$, where  ${\Cal  I}_X$  is  the  sheaf  of  ideals
of  ${\Cal  O}_{{\Bbb  P}^N}$ and Castelnuovo-Mumford regularity 
$\reg{(X)}$
of $X$  is defined as the least such $k$. 
There is a well-known conjecture concerning $k$-regularity: 
$ \reg{(X)}  \le\deg{(X)}-\codim{(X)}+1.$
This regularity conjecture including the classification of 
borderline examples was verified for integral curves (Castelnuovo,
Gruson, Lazarsfeld and Peskine), and an optimal bound was also
obtained for smooth surfaces (Pinkham, Lazarsfeld). It will be shown
here that $\reg{(X)}  \le\deg{(X)}-1$ for smooth threefolds $X$ in
$\Bbb P^5$ and that the only extremal cases are the Segre
threefold and the complete intersection of two quadrics. Furthermore,
every smooth threefold $X$ in ${\Bbb P}^5$ is $k$-normal \text{for
all} \,\,$k \ge \deg{(X)}-4$, which is the optimal bound as the
Palatini $3$-fold of degree $7$ shows. The same bound also holds for
smooth regular surfaces in $\Bbb P^4$ other than for the Veronese
surface.
\endabstract
\endtopmatter
\parindent=5mm
\document
\baselineskip=12pt

\subheading {0. Introduction}
 
Let  $X$  be  a projective  variety  of  dimension  $n$ and  degree
$d$  in  ${\Bbb P}^N$ defined over the complex number field ${\Bbb
C}$. $X$  is said to be $k$-normal  if  the  homomorphism $$ H^0
({\Bbb  P}^N , {\Cal  O}_{\Bbb  P}(k))  \rightarrow H^0 (X,  {\Cal
O}_X(k)) $$ is  surjective,  i.e.,  hypersurfaces  of  degree $k$  cut
out  a  complete linear system on $X$. We know that $X$ is $k$-normal
for all $k>>0.$ It would be very interesting to find an {\it explicit}
bound $k_0$ in terms of its invariants such that $X$ is $k$-normal for
all $k \ge k_0$ for all nondegenerate projective varieties.

On the other hand, as D.Mumford defined, we say that  $X$  is
$m$-regular  if   $H^i ({\Bbb  P}^N , {\Cal  I}_X(m-i))  =  0  \quad
\text{for  all}   \,\,  i  \ge  1$, where  ${\Cal  I}_X$  is  the
sheaf  of  ideals  of ${\Cal  O}_{{\Bbb  P}^N}$. It  is  easy  to
check  that  $X$  is  $(m+1)$-regular  if  and  only if  $X$  is
$m$-normal  and $H^i (X , {\Cal  O}_X(m-i))  =  0  \quad  \text{for
all }  i   > 0$.  Let's define   $\reg{(X)}$ to be  $ \min  \{ m  \in
{\Bbb   Z}  \:    X  \,\,  \text{is}  \,\,  m\text{-regular} \}$. More
generally,   for  a  coherent  sheaf  $\Cal F$ on ${\Bbb  P}^N$,
$\Cal F$  is  $m$-regular  if $\,\,  H^i ({\Bbb P}^N ,  {\Cal
F}(m-i))  =  0  \quad  \text{for all}  \,\,  i \ge 1$ $\,\,$  and
$\,\,  \reg{({\Cal  F})} \,\,$  is  defined  by  $\inf  \{ m   \in
{\Bbb  Z}  \:  {\Cal  F} \,\,  \text{is}  \,\,   m\text{-regular} \}$.

We  know  that if  $X$  is  $m$-regular,  then the degrees of the
minimal generators of the  saturated  ideal  $I_X$ defining $X$ are
bounded by $m$ and  hence  there  is no  $(m+1)$-secant line of
$X$.(More generally, the degrees of generators  of the $i$-th module
of syzygies  are bounded by $m+i$ \cite {EG}.)

Furthermore,  the  Hilbert polynomial  and  the  Hilbert  function  of
$X$  have  the  same values  for all $k  \geq  m-1.$
 
\proclaim{Regularity conjecture \cite{EG}, \cite{GLP}} \roster \item
$X$  is  $m$-normal $\text{for  all}  \,\,  m \ge d-\codim{(X)}$.
\item $X$ is $m$-regular $\text{for  all}  \,\,  m  \ge  d-\codim{(X)}
+ 1$, i.e., \newline $\reg{(X)}  \le  d- (N-n) + 1$. \endroster
\endproclaim In addition, the problem classifying all extremal
examples which make the bound best possible has been one of
interesting themes along these  lines.

Historically,  the  curve  case  in  $\Bbb P^3$  was  first  proved
by Castelnuovo  in 1893 \cite {C} and his result was completed by
Gruson, Lazarsfeld and Peskine in 1983 \cite{GLP}, where they proved ,
by  using  vector bundle  techniques, that every integral curve $X$ of
degree $d$  in $\Bbb P^N$ is $(d-\codim{(X)} + 1)$-regular. They also
made  the  list of all  extremal curves  which  fail  to  be
$(d-\codim{(X)})$-regular
 
For smooth surfaces, bounds are obtained in  \cite{P} and \cite{L},
where  R. Lazarsfeld got the optimal bound $\reg{(X)}  \le
(d-\codim{(X)} +1)$. For smooth threefolds and smooth fourfolds,
bounds are obtained in \cite{K}, where the author got the following
results: $\reg{(X)} \le (d-\codim{(X)} + 2)$ for smooth threefolds and
$\reg{(X)} \le (d-\codim{(X)}+5)$ for smooth fourfolds. Note that
these bounds are  off by 1 (resp., 4) from the conjectured bound. For
arbitrary  dimensional  smooth  projective  varieties, $ \reg{(X)}
\le  \min \{e,n+1 \} (d-1) - n+1, e=\codim{(X)} $ by Bertram, Ein and
Lazarsfeld  \cite{BEL}. In the case of arbitrary codimension two
smooth projective variety $X$ of degree  $d$  in  ${\Bbb P}^N$, Alzati
\cite{Al} got a better bound, i.e., $ \reg{(X)}  \le (d+2)+
\frac{1}{2}N(N-1)-2N $ (especially, $\reg{(X)}  \le d+
\frac{1}{2}N(N-1)-(N-1)[\frac{(N+4)}{2}],  N \ge 6 )$ by using various
facts on higher order normality and vanishing theorems for an ideal
sheaf due to Alzati and Ottaviani and recently, Peeva and Sturmfels
showed in  \cite{PSt} that $ \reg{(X)}  \le \deg{(X)} $ for toric
varieties of codimension two. However, unfortunately, for smooth
projective varieties of dimension $n \ge 2$, extremal examples with
geometric interpretations has not  been classified yet.

The purpose of this paper is to give new bounds for regularity of
smooth threefolds in  $\Bbb P^5$ and to classify extremal examples by
using local cohomology.

\proclaim{Theorem 1} If a locally Cohen-Macaulay projective variety
$X$ of dimension $n$ is contained in a hyperquadric $Q$ in $\Bbb
P^{n+2}$ then $X$  is either  a  complete  intersection  or
projectively  Cohen-Macaulay and  linked  to  a  linear  space $L
\simeq  {\Bbb  P}^n$ via  a  complete  intersection  of $Q$  and
another  hypersurface. \endproclaim \proclaim{Theorem 2} Let $X$ be a
smooth threefold of degree $d$ in $\Bbb P^5$. Then $\reg{(X)} \le
d-1$. Furthermore, the only extremal examples are the Segre threefold
and  the complete intersection of two quadrics.
\endproclaim

\proclaim{Corollary 3} Any smooth threefold in $\Bbb P^5$ is
$k$-normal for all $k \ge d-4$. The Palatini scroll of degree $7$
which is not 2-normal(see example 2.0.) shows that this is the optimal
bound. In addition, the same is true for regular surfaces in $\Bbb
P^4$ except for the Veronese surface. \endproclaim

{\bf Acknowledgements.} I would like to express my gratitude to
Professor Henry Pinkham for his encouragement and stimulating
discussion with geometric insight throughout my graduate studies at
Columbia University. In addition, I also thank my colleagues,
especially Bo Ilic,  for valuable discussions in the seminar group
during my stay there.

\subhead {1. Basic background}\endsubhead
In  this  section,  we  recall  the  definitions  and  well-known
results  which are  used  in  the  following  sections.

\proclaim{Lemma 1.1}
Let $\Cal F$  be  a  $p$-regular  vector  bundle  on ${\Bbb P}^n$  and
$\Cal G$  be  a  $q$-regular  vector  bundle  on ${\Bbb  P}^n$ defined
over the complex numbers  ${\Bbb C}$. Then ${\Cal  F}  \otimes  {\Cal
G}$  is  $(p+q)$-regular  and  $S^k({\Cal F})$ and $\bigwedge^k({\Cal
F})$  are  $(kp)$-regular.
\endproclaim

\demo{Proof} See [L], p.428. \enddemo

Let  $X$  be  a  locally  Noetherian  scheme, $Y$ a  closed  subscheme
of  $X$, and  let  $\Cal   F$  be  a  coherent  sheaf  on  $X$. We
define $$ \depth_Y  {\Cal  F}  =  \inf_{x\in Y}  (\depth{{\Cal
F}_x}), $$ where  $x$  is  not  necessarily  a  closed  point.

Consider  the  following  left  exact  functor $\underline \Gamma_Y(X,
\cdot)$ from  the  category  of  sheaves  of  abelian groups to  the
category  of sheaves  of  abelian  groups,  defined  by $$ \underline
\Gamma_Y(X, \Cal F) (U)=\Gamma_{Y\cap U}(U,\Cal F\mid_U) $$ which  is
a  subgroup  of  sections  of  $\Cal F\mid_U$ with  support  in
$Y\cap U$. Then  the  $p$-th  right  derived  functors  of $\underline
\Gamma_Y(X, \cdot)$ are  denoted  by $\underline H_Y^p(X,  \cdot),
p=0,1,\cdots$. Similarly, consider  the left  exact  functor
$\Gamma_Y(X, \cdot)$ from  the  category  of  sheaves  of  abelian
groups to  the  category  of  abelian  groups,  defined  by $\Cal F
\longmapsto  \Gamma_Y(X, \Cal F)$ which  is  a  subgroup of global
sections  of  $\Cal F$ with  support  in   $Y$. We denote the $p$-th
right  derived  functors  of $\Gamma_Y (X, \cdot)$ by $H_Y^p (X,
\cdot).$

Note that $\underline H_Y^p(X,  \cdot)$ is the sheaf associated to the
presheaf $$ U \longmapsto H_{Y\cap U}^p(U,\Cal F\mid_U) $$ and there
is a spectral sequence $H_Y^n(X,  \Cal F) \Longleftarrow E_2^{p,q}
=H^p(X, \underline H_Y^q(\Cal F)).$

\proclaim{Theorem 1.2}  
Let  $X$  be  a  locally  Noetherian  scheme, $Y$  a  closed
subscheme  of  $X$,   and let  $\Cal F$  be  a  coherent  sheaf  on
$X$. Then  the  following  conditions  are  equivalent: \roster
\item"(1)"  $\underline H_Y^i(X,\Cal F)=0$  for all $i<n$, \item"(2)"
$\depth_Y{\Cal F}\ge n$. \endroster \endproclaim

\demo{Proof} See \cite{Gro} p.44.
\enddemo

We  assume that  $X$  is  a  locally  Cohen-Macaulay subscheme  in
${\Bbb  P}^N$ over the field $\Bbb C$  and $ (M^i)(X)  =
H^i_{\ast}({\Cal  I}_X)  =\allowmathbreak \bigoplus_{m  \in  {\Bbb
Z}}  H^i({\Bbb  P}^N,  {\Cal  I}_X(m)) $ for $1  \leq  i  \leq
\dim{(X)}$ and $ [(M^i)(X)]^{\Bbb C}  = \Hom_{\Bbb C}{((M^i)(X), {\Bbb
C})} $. \proclaim{Theorem  1.3} Let  $X_1, X_2  \subset  {\Bbb  P}^N$
have  dimension  $n.$ Assume  that  $X_1$  is  linked  to  $X_2$  via
a  complete  intersection $X  =  Z(F_1,F_2, \dots, F_{N-n})$  with
$\deg{(F_i)}  =  d_i$. Let  $d  =  \sum  d_i$.  Then \roster \item[1]
$(M^{n-i+1})(X_2)  \simeq  [(M^i)(X_1)]^{\Bbb C}(N+1-d)$ \item[2] We
have  an  exact  sequence $$ 0  \rightarrow {\Cal  I}_X  \rightarrow
{\Cal  I}_{X_1}  \rightarrow \omega_{X_2}(N+1-d)  \rightarrow 0 $$
where  $\omega_{X_2}$  is  the  dualizing  sheaf  of  $X_2$.
\endroster \endproclaim

\demo{Proof}
See  \cite{Mi}, Ch IV. \enddemo

\proclaim{Theorem  1.4 (The  (dimension+2)-secant  Lemma)} Let  $X
\subset  {\Bbb  P}^N$  be  a  smooth  $n$-dimensional subvariety  and
let   $Y$  be  an  irreducible  variety parametrizing  a  family
$\{L_y\}$  of   lines  in   ${\Bbb  P}^N$. Assume  that,  for  a
general  $L_y$,  the  length  of  the scheme-theoretic  intersection
$L_y  \cap  X$  is  at  least $(n+2)$.  Then  we  have $$
\dim{(\cup_{y  \in  Y}  L_y)}  \leq  n+1 $$ \endproclaim

\demo{Proof}
See  \cite{Ran2}.
\enddemo

\subhead  {2. Castelnuovo-Mumford Regularity for Smooth Threefolds  in
${\Bbb  P}^5$ and Extremal Examples} \endsubhead

Throughout this section, we work over the complex number field  ${\Bbb
C}$. A conjecture of R.~Hartshorne (which is still open)  states  that
any smooth $n$-manifold $X$  in  ${\Bbb  P}^{n+2}, (n \  \geq  \  4)$
is  a  complete intersection, equivalently  it is projectively normal.
However,  we  have  many  examples  of  smooth non-projectively normal
surfaces  in  ${\Bbb  P}^4$  and  3-folds  in  ${\Bbb P}^5$,   some of
which  are  defined  scheme-theoretically  as  the  dependency loci of
a morphism   $\varphi : {\Cal  F} \rightarrow  {\Cal  G}$ between
vector  bundles  ${\Cal  F}$  and  ${\Cal  G}$  of rank $m$  and $m+1$
on  ${\Bbb  P}^N$,  respectively.  In  other words, let  $\SHom{(\Cal
F,\Cal  G)}$  be  globally generated  and  $\varphi$  be  a generic
morphism.  Then  the morphism  $\varphi$  induces
$$
(\bigwedge^m ({\Cal  G}))^{\vee} \otimes (\det{\Cal  G}) \simeq
{\Cal  G}  @> (\bigwedge^m \varphi)^{\ast} >>
(\bigwedge^m ({\Cal  F}))^{\vee} \otimes (\det{\Cal  G})  \simeq
{\Cal  O}_{{\Bbb  P}^N}(c_1({\Cal  G}) - c_1({\Cal  F}))
$$
which  fits  into  an  exact  sequence $$ 0 \rightarrow {\Cal  F} @>
\varphi >> {\Cal  G} \rightarrow {\Cal  I}_X(c_1({\Cal  G}) -
c_1({\Cal  F})) \rightarrow 0, $$ where  X  is  the  dependency  locus
of  $\varphi : {\Cal  F} \rightarrow {\Cal  G}$.  By   Kleiman's
Bertini-type  theorem \cite{Klm}, $X$  is  a  codimension two
subvariety in  ${\Bbb  P}^n$  and nonsingular away from a subset of
codimension $\le 4 $ in $X$.
 
\subhead{Examples 2.0}
\endsubhead 
From the exact sequence $$ 0 \rightarrow \Omega^2_{{\Bbb  P}^5}(2+t)
\rightarrow \bigwedge^2 ({\italic  V})\otimes {\Cal  O}_{{\Bbb  P}^5}(t)
\rightarrow \Omega_{{\Bbb  P}^5}(2+t)  \rightarrow 0, $$ we know that
$\Omega_{{\Bbb  P}^5}(2+t)$  is  globally  generated for all $t \ge
0$. Then,  four generic  sections  $s_1,  s_2, s_3,  s_4$ of
$H^0({\Bbb  P}^5,  \Omega_{{\Bbb  P}^5}(2+t))$ induce  an  exact
sequence $$ 0  \rightarrow {\Cal  O}_{{\Bbb  P}^5}^{\oplus 4}  @>
\varphi=(s_1,s_2,s_3,s_4) >> {\Omega}_{{\Bbb  P}^5}^1(2+t)  @>
(\bigwedge^4 \varphi)^{\ast} >> {\Cal  I}_{X_t}(c_1({\Omega}_{{\Bbb
P}^5}(2+t)))  \rightarrow 0 $$ where $X_t$  is  the  dependency  locus
of  $\varphi$.  $X_t$ is smooth by Kleiman's Bertini-type theorem
\cite{Klm}.        By the formula \cite{OSS} p.16, $$ c_k(E \otimes
{\Cal  L}) = \sum_{i=0}^{k} \binom {r-i}{k-i} c_i(E) \cdot c_1({\Cal
L})^{k-i},\,  r =\rank{(E)},\,\, {\Cal  L}\,\,\text
{a\,\,line\,\,bundle}, $$ we  get $ c_1(\Omega_{{\Bbb  P}^5}^1(2+t)) =
4 + 5t, $ $ c_2(\Omega_{{\Bbb  P}^5}^1(2+t)) = 10t^2 + 16t + 7 =
\deg{(X_t)}, $ and we can compute $$ h^1({\Bbb  P}^5,  {\Cal
I}_{X_t}(k))  = \left\{ \aligned &1, \  when \  k = 4t + 2   \\ &0, \
when \  k \neq 4t + 2 \endaligned \right. $$ via the vanishing of
$h^q({\Bbb  P}^n,\Omega^p_{{\Bbb  P}^n}(k))$ \cite{OSS}. So,  $X_t$
is  $k$-normal  for  all  $k  \neq  4t + 2$. Note that $X_0$ is called
the Palatini scroll. It has  degree $7$ and sectional genus 4.
According to a conjecture of Peskine and Van de Ven, it is the unique,
non 2-normal, smooth threefold in $\Bbb P^5.$

\proclaim{Theorem  2.1}
If   a  locally  Cohen-Macaulay,  nondegenerate  subvariety $X^n
\subset \Bbb P^{n+2}(n \geq 2)$ is  contained  in  a hyperquadric $Q$,
then  $X$  is  either  a  complete  intersection  or is projectively
Cohen-Macaulay and  linked  to  a  linear  space  $L  \simeq  {\Bbb
P}^n$ via  a  complete  intersection  of $ Q$  and   another
hypersurface. \endproclaim

\demo{Proof}
By  a  suitable  change  of  coordinates,  we  may  assume  that  $Q$
is  defined  by  a  homogeneous  polynomial $x_0^2 + x_1^2 + \dots +
x_k^2$  of  rank  $k+1$. Let's consider  all  cases  according  to
$\rank{(Q)}$. \roster \item"Case {\rm  I}: " $\rank{(Q)}   =   k+1
\geq  5$. In this case, $X$ is a complete intersection. See \cite{H2},
exercise II.6.5(d), Klein's theorem. \item"Case  {\rm  II}: "
$\rank{(Q)}   =  4$. Then  $Q$  is  the  cone  over  a  smooth
quadric  surface  $Q_1$  in ${\Bbb  P}^3$  with  vertex  $\Lambda
\cong {\Bbb  P}^{n-2}$. Let ${\Cal U}={Q - \Lambda}$. Then  $X -
\Lambda$  is  a  Cartier  divisor  of  $\Cal  U$, so ${\Cal  I}_{X /
Q}|_{\Cal  U}$  is  a  line  bundle  on ${\Cal  U}$. We have an
isomorphism ${\pi}^{\ast}_{\Lambda} : Cl(Q_1) @> \sim >> Cl({\Cal
U})=Cl(Q)$ where the projection  ${\pi}_{\Lambda} : {\Cal  U}
\rightarrow  Q_1 \subset  {\Bbb  P}^3$ is  an  affine  bundle with
fibre  ${\Bbb  A}^{n-1}$ and $\Pic{(Q_1)}=\Pic{({\Cal U})} \simeq
{\Bbb  Z} \oplus {\Bbb  Z}$.  We'll show that ${\Cal  I}_{X /
Q}|_{\Cal  U}$ is of type (a,b), $|a-b|  \leq  1$,  where ${\Cal
I}_{X  /  Q}|_{\Cal  U}  =  {\pi}_{\Lambda}^{\ast}({\Cal   O}_{Q_1}
(a,b))$ for some $ a,b \in \Bbb Z.$ The proof of case II will be
deferred until after some auxiliary propositions and lemmas.
\endroster

\proclaim{Lemma 2.2} 
Let $X^n  \subset \Bbb  P^{n+2},\,  (n \geq 2)$  be  a  locally
Cohen-Macaulay subvariety contained  in  a hyperquadric  $Q$ of rank
$4.$ Then  $\depth_{\Lambda}{( {\Cal  I}_{X / Q} )} \geq 3$.
\endproclaim
\demo{Proof} Since ${\Lambda}^{n-2} \subset Q \subset
{\Bbb  P}^{n+2}$ and  $Q$  is  locally  Cohen-Macaulay,
$$
\depth_{\Lambda}{({\Cal  O}_Q)}  = \inf_{x \in \Lambda}{\depth{( {\Cal
O}_{Q,x} )}} = \depth{( {\Cal  O}_{Q, \eta} )} = \height{\eta} = 3,
$$
where  $\eta$  is  the  generic  point  of  $\Lambda^{n-2}.$ Since
$n-4 \leq \dim{Z}\leq n-2$ for any irreducible component $Z$ of $X
\cap \Lambda $, we get $2 \le\depth{({\Cal  O}_{X,\wp})}=\dim{ ({\Cal
O}_{X,\wp})} \leq  4 $ for any minimal prime ideal  $\wp$  containing
$I_X$  and $I_{\Lambda}$. Therefore,  $\depth_{\Lambda}{( {\Cal  O}_X
)}  \geq  2 $  whenever $X \cap \Lambda \neq \phi$. From a  long
exact  sequence  of  local  cohomology  groups
$$
\multline
\underline{H}^0_{\Lambda}(Q,{\Cal  I}_{X /
Q})  @>>> \underline{H}^0_{\Lambda}(Q,{\Cal  O}_Q)  @>>>
\underline{H}^0_{\Lambda}(Q,{\Cal  O}_X)\\  @> \delta >>
\underline{H}^1_{\Lambda}(Q,{\Cal  I}_{X / Q})  @>>>
\underline{H}^1_{\Lambda}(Q,{\Cal  O}_Q)  @>>>
\underline{H}^1_{\Lambda}(Q,{\Cal  O}_X)\\  @>>>
\underline{H}^2_{\Lambda}(Q,{\Cal  I}_{X / Q}) @>>>
\underline{H}^2_{\Lambda}(Q,{\Cal  O}_Q)  @>>>
\underline{H}^2_{\Lambda}(Q,{\Cal  O}_X) @>>>
\endmultline
\tag{$2.2.1$}
$$ Since  $\depth_{\Lambda}{( {\Cal  O}_Q )} = 3$  and
$\depth_{\Lambda}{( {\Cal O}_X )} \geq  2$, by using Theorem  1.2.,
$$\depth_{\Lambda}{({\Cal  I}_{X / Q})} \geq\min{(\depth_{\Lambda}{(
{\Cal  O}_Q )}, \depth_{\Lambda}{( {\Cal O}_X )}+1 )}=3.$$ \qed
\enddemo

\proclaim{Corollary 2.3}  We  have  isomorphisms  $H^i(Q,{\Cal  I}_{X
/ Q}(m)) \simeq  H^i({\Cal  U},{\Cal  I}_{X / Q}(m)|_{\Cal  U}),$  for
all $m \in \Bbb Z \,$ and $\ i = 0,1$. \endproclaim

\demo{Proof}
For  a  coherent  sheaf ${\Cal  I}_{X / Q}(m)$  on  $Q$,  we  have  a
long  exact  sequence  of cohomology  groups $$ \aligned 0 @>>> \\
@>>>\\@>>> \\ \endaligned \aligned &  H^0_{\Lambda}(Q,{\Cal  I}_{X /
Q}(m)) \\ &  H^1_{\Lambda}(Q,{\Cal  I}_{X / Q}(m)) \\ &
H^2_{\Lambda}(Q,{\Cal  I}_{X / Q}(m)) \\ \endaligned \aligned &  @>>>
H^0(Q,{\Cal  I}_{X / Q}(m))  \\ &  @>>>  H^1(Q,{\Cal  I}_{X / Q}(m))
\\ &  @>>>  H^2(Q,{\Cal  I}_{X / Q}(m))  \\ \endaligned \aligned &
@>>>  H^0({\Cal  U},{\Cal  I}_{X / Q}(m)|_{\Cal  U})  \\ &  @>>>
H^1({\Cal  U},{\Cal  I}_{X / Q}(m)|_{\Cal  U})  \\ &  @>>>  H^2({\Cal
U},{\Cal  I}_{X / Q}(m)|_{\Cal  U})  \\ \endaligned \tag{$2.3.1$} $$
By Theorem 1.2 and Lemma  2.2,   $\underline{H}^i_{\Lambda}(Q,{\Cal
I}_{X / Q}(m)) = 0$ for $i = 0,1,2.$ Therefore,  we   have
isomorphisms  $H^i(Q,{\Cal  I}_{X / Q}(m)) \simeq  H^i({\Cal  U},{\Cal
I}_{X / Q}(m)|_{\Cal  U}),  \,\,  i = 0,1$. \qed
\enddemo

\proclaim{Lemma  2.4}
Let $Q = \Proj{{\Bbb  C}[x_0,x_1,\dots,x_{n+2}]} / (x_0 x_3 - x_1 x_2)$
be  a  hyperquadric   of rank $ 4$  in  ${\Bbb  P}^{n+2}$ with  a
vertex  ${\Lambda}  =  Z(x_0,x_1,x_2,x_3)$  and ${\Lambda} = {\Bbb
P}(V),  V =  {\Bbb C}x_4 \oplus  {\Bbb C}x_5 \oplus \dots \oplus
{\Bbb C}x_{n+2}$. Let  ${\pi}_{\Lambda} : {\Cal  U} = Q - \Lambda
\rightarrow Q_1$ be a  projection from  the  center  $\Lambda$.  Then
$ {\pi}_{{\Lambda}_{\ast}} {\Cal  O}_{\Cal  U} = {\bigoplus}_{m =
0}^{\infty} S^m V \otimes {\Cal  O}_{Q_1}(-m) $. \endproclaim

\demo{Proof}
Let  $\{ D^+(x_i) | i = 0,1,2,3\}$  be  the  standard  covering  of
$Q_1$ and $A_i = {\Bbb  C}[\frac{x_0}{x_i}, \dots,\frac{x_3}{x_i}] /
(\frac{x_0}{x_i} \frac{x_3}{x_i} - \frac{x_1}{x_i} \frac{x_2}{x_i} )$
be  the  affine  coordinate  ring  of  $Q_1 \cap D^+(x_i)$ for $i  =
0,1,2,3$. Then $$ {\pi}_{\Lambda}|_{D^+(x_i)} : Spec A_i
[\frac{x_4}{x_i},\frac{x_5}{x_i},\dots,\frac{x_{n+2}}{x_i}]
\rightarrow Spec A_i $$ is  a  restriction  of  ${\pi}_{\Lambda}$. On
$D^+(x_i) \cap D^+(x_j)$, $ (\frac{x_k}{x_i})^m = (\frac{x_j}{x_i})^m
\cdot (\frac{x_k}{x_j})^m, 0 \leq i,j \leq 3, and \,\, 4 \leq k \leq
n+2. $ Furthermore,  $(\frac{x_j}{x_i})^m$  is  the  transition
function  of ${\Cal  O}_{Q_1}(-m)$.  Therefore, $$
{\pi}_{{\Lambda}_{\ast}} {\Cal  O}_{\Cal  U}  =
{\bigoplus}_{m=0}^{\infty} S^m(V) \otimes  {\Cal  O}_{Q_1}(-m), $$
where  $\Lambda = {\Bbb  P}(V)$  and  $S^m(V)$  is  the  $m$-th
symmetric  power of  $V =  {\Bbb C}x_4 \oplus \dots \oplus  {\Bbb
C}x_{n+2}$. This completes the  proof  of  Lemma  2.4. \qed
\enddemo

\proclaim{Proposition  2.5} Under the same situation as Lemma 2.2,
\roster
\item[1]  
$|a-b|  \leq  1$,  where  
${\Cal  I}_{X  /  Q}|_{\Cal  U}  \simeq 
{\pi}_{\Lambda}^{\ast}({\Cal   O}_{Q_1}
(a,b))$
\item[2]  
$H^1(Q,{\Cal  I}_{X / Q}(k)) = 0  {\text  \   for  \   all}  \  k \in  {\Bbb  Z}$
\endroster
\endproclaim

\demo{Proof}
Under the isomorphisms ${\pi}^{\ast}_{\Lambda} : Cl(Q_1) @> \sim >>
Cl({\Cal  U}) =  Cl(Q)$ and $Cl(Q_1) \simeq  Cl({\Bbb  P}^1  \times
{\Bbb  P}),$  we have ${\Cal  I}_{X / Q}|_{\Cal  U} =
{\pi}^{\ast}_{\Lambda}({\Cal O}_{Q_1}(a,b))$ for  some  $a,b \in {\Bbb
Z}$  and ${\Cal  I}_{X / Q}(k)|_{\Cal  U} =
{\pi}^{\ast}_{\Lambda}({\Cal  O}_ {Q_1}(a+k,b+k))$. Since
${\pi}_{\Lambda}$  is  an  affine  morphism, so $R^i
{\pi}_{{\Lambda}_{\ast}} ({\Cal O}_{\Cal  U}) = 0,  \,\,  i > 0  \,\,$
and $H^i({\Cal  U},{\Cal  I}_{X / Q}(k)|_{\Cal  U})  = H^i({\Cal
U},{\pi}^{\ast}_{\Lambda} {\Cal  O}_{Q_1}(a+k,b+k))  = H^i(Q_1,{\Cal
O}_{Q_1}(a+k,b+k) \otimes {\pi}_{{\Lambda}_{\ast}} {\Cal  O}_{\Cal
U})$ for  all  $i  \geq  0$. Therefore,  for  all  $i = 0, 1$, $$
\aligned &\quad\,\,\,\,H^i(Q,  {\Cal  I}_{X  /  Q}(k))  \cong
H^i({\Cal  U},  {\Cal  I}_{X  /  Q}(k)|_{\Cal  U})  \\ &\cong
H^i(Q_1,{\Cal  O}_{Q_1}(a+k,b+k) \otimes {\pi}_{{\Lambda}_{\ast}}
{\Cal O}_{\Cal  U})  \\ &\cong  {\bigoplus}_{m=0}^{\infty} S^m(V)
\otimes H^i(Q_1,{\Cal O}_{Q_1}(a+k-m,b+k-m)). \endaligned $$ Since
$Q_{1} \simeq  {\Bbb  P}^1  \times  {\Bbb  P}^1$, it follows
immediately from the K{\" u}nneth formula, $$ \spreadlines{12pt}
\aligned &\quad\,\,\,\,  H^1(Q,  {\Cal  I}_{X  /   Q}(k))  \\ &\cong
{\bigoplus}_{m=0}^{\infty}  S^m(V)  \otimes H^1(Q_1,{\Cal
O}_{Q_1}(a+k-m,b+k-m))  \\ &\cong  {\bigoplus}_{m=0}^{\infty} \{S^m(V)
\otimes H^1({\Bbb  P}^1,{\Cal  O}_{{\Bbb  P}^1}(a+k-m) \otimes
H^0({\Bbb  P}^1,{\Cal  O}_{{\Bbb  P}^1}(b+k-m))\} \oplus  \\
&\qquad\qquad\,\,\,\,\{S^m(V)  \otimes H^0({\Bbb  P}^1,{\Cal
O}_{{\Bbb  P}^1}(a+k-m)  \otimes H^1({\Bbb  P}^1,{\Cal  O}_{{\Bbb
P}^1}(b+k-m)) \} \endaligned $$ By  Serre's  vanishing  Theorem,
$H^1(Q,{\Cal  I}_{X / Q}(k)) = 0$  for $k >> 0$.  Therefore, $$
H^1(Q_1,{\Cal O}_{Q_1}(a+k-m,b+k-m)) = 0, \,\,  \text{for  all}  \,\,
m  \in  {\Bbb  Z}^+   \,\,   \text{and}  \,\,  \text{for  all}  \,\,
k >> 0. $$ Let's fix an integer $k_0$  such  that  $H^1(Q,{\Cal  I}_{X
/ Q}(k_0)) = 0$ and assume $ b \ge a$. We  know  that in  general, $$
\aligned &h^0({\Cal  O}_{{\Bbb  P}^1}(b+k_0-m)) = h^0({\Cal  O}_{{\Bbb
P}^1}(a+k_0-m)) = 0, \,\,\text{for  all}  \,\,   m  \geq  b+ k_0 + 1
\\ &h^1({\Cal  O}_{{\Bbb  P}^1}(b+k_0-m)) = h^1({\Cal  O}_{{\Bbb
P}^1}(a+k_0-m)) = 0, \,\,\text{for  all}  \,\,  m  \leq  a + k_0 + 1
\endaligned $$ If  $a+2 \le b$,  then $\{m | a+k_0+2 \leq  m \leq
b+k_0 \}$  is  not empty.  This  implies  that $$ H^1(Q_1,{\Cal
O}_{Q_1}(a+k_0-m,b+k_0-m)) \neq 0  \,\,  for  \,\,  m = a + k_0 +2, $$
which  contradicts  the  fact  that $H^1(Q,{\Cal  I}_{X / Q}(k_0))
=0$.  So $|a-b| \leq 1$. In  particular,  by  repeating  this
argument, $H^1(Q,{\Cal  I}_{X / Q}(k))  = 0$  for  all  $k \in  {\Bbb
Z}$. \qed
\enddemo

Now  let's  go  back  to  the  proof  of  Case  {\rm  II} in Theorem
2.1. Consider the exact sequence $$ 0 \rightarrow {\Cal  I}_{Q / {\Bbb
P}^{n+2}}(k) \rightarrow {\Cal  I}_{X / {\Bbb  P}^{n+2}}(k)
\rightarrow {\Cal  I}_{X / Q}(k) \rightarrow 0 $$ Since  $h^1({\Bbb
P}^{n+2}_{\Bbb C},  {\Cal  I}_{Q / {\Bbb P}^{n+2}}(k)) = h^1(Q,{\Cal
I}_{X  /  Q}(k))  =  0  $  for  all  $  k \in  {\Bbb Z}$, we  have
$h^1({\Cal  I}_{X / {\Bbb  P}^{n+2}}(k))  = 0$     for  all $   k  \in
{\Bbb  Z}$. Furthermore,  if $ {\Cal  I}_{X  /  Q}|_{\Cal  U}  \cong
{\pi}_{\Lambda}^{\ast}{\Cal  O}_{Q_1}(-a,-a)$ for some $a > 0, $ then
$X - {\Lambda}$  is  linearly equivalent  to  $(a, a)$-type  in
${\italic  Cl}({\Cal  U})  \cong  {\Bbb  Z}  \oplus  {\Bbb  Z}$.  So
$X$  is also  linearly  equivalent  to  a class of type $(a,a)$  in
${\italic  Cl}(Q)$  and $\deg{(X)}  =  2a$.  On  the  other  hand,
we  know  that $$ \aligned &\quad\,\,\,\,H^0(Q,  {\Cal  I}_{X  /
Q}(a))  \cong H^0({\Cal  U},  {\Cal  I}_{X  /  Q}(a)|_{\Cal  U})  \\
&\cong  \bigoplus_{m=0}^{\infty} S^m(V)  \otimes H^0(Q_1,  {\Cal
O}_{Q_1}(-m,-m))  \\ &\cong  \bigoplus_{m=0}^{\infty} S^m(V)  \otimes
H^0({\Bbb  P}^1,  {\Cal  O}_{{\Bbb  P}^1}(-m))  \otimes H^0({\Bbb
P}^1,  {\Cal  O}_{{\Bbb  P}^1}(-m))  \\ &\cong\,\,\, {\Bbb C}
\endaligned $$
From  the exact  sequence
$$ 0  \to H^0({\Bbb P}^{n+2},  {\Cal  O}_{{\Bbb  P}^{n+2}}(a-2))
\overset{\times  Q}\to\longrightarrow
H^0({\Bbb  P}^{n+2},  {\Cal  I}_{X  /  {\Bbb  P}^{n+2}}(a)) \to H^0(Q,
{\Cal  I}_{X  /  Q}(a))  \to 0 $$ we  get $ h^0({\Bbb  P}^{n+2},
{\Cal  I}_{X  /  {\Bbb  P}^{n+2}}(a))  = h^0({\Bbb  P}^{n+2},  {\Cal
O}_{{\Bbb  P}^{n+2}}(a-2))  + 1 $ and  $X$  is  a  complete
intersection  of  $Q$  and another  hypersurface  of  degree $a$.

Similarly,  if ${\Cal  I}_{X  /  Q}|_{\Cal  U}  \cong
{\pi}_{{\Lambda}^{\ast}} {\Cal  O}_{Q_1}(-a, -a+1)$ in $\Pic{(\Cal
U)}$, then  $X$  is  linearly  equivalent  to  $(a, a-1)$-type  in
${\italic  Cl}(Q)$  and  $\deg{(X)}  =  2a-1$. Note that $$ \aligned
&\quad\,\,\,\,H^0(Q,  {\Cal  I}_{X  /  Q}(a))  \cong H^0({\Cal  U},
{\Cal  I}_{X  /  Q}(a)|_{\Cal  U})  \\ &\cong
\bigoplus_{m=0}^{\infty} S^m(V)  \otimes H^0(Q_1,  {\Cal
O}_{Q_1}(-m,-m+1))  \\ &\cong  \bigoplus_{m=0}^{\infty} S^m(V)
\otimes H^0({\Bbb  P}^1,  {\Cal  O}_{{\Bbb  P}^1}(-m))  \otimes
H^0({\Bbb  P}^1,  {\Cal  O}_{{\Bbb  P}^1}(-m+1))  \\ &\cong   {\Bbb C}
\oplus   {\Bbb C} \endaligned \tag{$\ast$} $$ and $ h^0({\Bbb
P}^{n+2},  {\Cal  I}_{X  /  {\Bbb  P}^{n+2}}(a))  = h^0({\Bbb
P}^{n+2},  {\Cal  O}_{{\Bbb  P}^{n+2}}(a-2))  + 2. $ Hence  $X$  is
linked  to  $L  \cong  {\Bbb  P}^n$  via  a  complete intersection  of
$Q$  and  another  hypersurface  of  degree  $a$. By  the  way,  the
property  of  being  projectively  Cohen-Macaulay is  invariant  in  a
liaison  class  by  Theorem  1.3,(1). Thus  $X$  is projectively
Cohen-Macaulay  because  $L^n$  is  projectively  Cohen-Macaulay.

Case  {\rm  III}:  $\rank{(Q)}  =  3$.  We may assume that $V  = {\Bbb
C} x_3  \oplus  {\Bbb C}x_4  \oplus  \dots  \oplus  {\Bbb C}x_{n+2}$,
$\Sing{(Q)}=\Lambda^{n-1}  = {\Bbb  P}(V)$, $Q = <C, \Lambda>$, where
$C$  is a  smooth  conic  in  ${\Bbb  P}^2$. Since $X$ is locally
Cohen-Macaulay,  we can get $\depth_{\Lambda}{({\Cal O}_Q)}  =  2$ and
$\depth_{\Lambda}{{\Cal  O}_X}  \ge  1$. Hence, by ($2.2.1$), $
\underline{H}^i_{\Lambda}(Q,  {\Cal  I}_{X  /  Q})  =  0 \,\,
\text{for}  \,\,  i  =  0,1 $ and $ H^0(Q,  {\Cal  I}_{X  /  Q}(m))
\cong H^0({\Cal  U},  {\Cal  I}_{X  /  Q}(m)|_{\Cal  U}) $. Note  that
$ {\italic  Cl}(Q)  \cong {\italic  Cl}({\Cal  U})  \cong {\italic
Cl}(C)  \cong {\italic  Cl}({\Bbb  P}^1)  \cong {\Bbb  Z} $.
If $\,\, {\Cal  I}_{X  /  Q}|_{\Cal  U}  =
{\pi}_{\Lambda}^{\ast}({\Cal  O}_C(-a))  = {\pi}_{\Lambda}^{\ast}(
i^{\ast}({\Cal  O}_{{\Bbb  P}^1}(-2a))) \,\,  \text{for  some}  \,\,
a > 0 $, where  $i : C  \rightarrow {\Bbb  P}^1$ is  an  isomorphism,
we have $\deg{(X)}=2a$ and
$$
H^0(Q,{\Cal  I}_{X  /  Q}(a))  \cong
H^0({\Cal  U},  {\Cal  I}_{X/Q}(a)|_{\Cal  U}) 
\cong
\bigoplus_{m=0}^{\infty}  S^m(V)
\otimes  H^0(C,  {\Cal  O}_C(-m))
\cong{\Bbb C},
$$
and  the  exact  sequence
$$
0  \rightarrow H^0({\Bbb P}^{n+2},  {\Cal  O}_{{\Bbb  P}^{n+2}}(a-2))
\overset{\times Q}\to\longrightarrow H^0({\Bbb  P}^{n+2},
{\Cal  I}_{X  /  {\Bbb  P}^{n+2}}(a))
\rightarrow H^0(Q,  {\Cal  I}_{X  /  Q}(a))  \rightarrow 0
$$
implies
that $X$  is  a  complete   intersection  of  $Q$  and another
hypersurface F of  degree  $a$.

Next,  if $\,\, {\Cal  I}_{X  /  Q}|_{\Cal  U}  \cong
{\pi}_{\Lambda}^{\ast}(  i^{\ast}{\Cal  O}_{{\Bbb  P}^1}(-2a+1)) $,
then  $\deg{(X)}  =  2a  -1$  and $$ \spreadlines{12pt} \aligned
&\quad\quad H^0(Q,  {\Cal  I}_{X  /  Q}(a))  \cong H^0({\Cal  U},
{\Cal  I}_{X  /  Q}(a)|_{\Cal  U})\\ &\cong H^0(C, i^{\ast}{\Cal
O}_{{\Bbb  P}^1}(-2a+1) \otimes {\Cal  O}_{C}(a) \otimes
{\pi}_{{\Lambda}_{\ast}} {\Cal O}_{\Cal  U}) \\ &\cong
\bigoplus_{m=0}^{\infty}  S^m(V)  \otimes H^0(C,  i^{\ast}{\Cal
O}_{{\Bbb  P}^1}(-2a+1) \otimes {\Cal  O}_{C}(a-m))\\ &\cong
\bigoplus_{m=0}^{\infty}  S^m(V)  \otimes H^0({\Bbb  P}^1,  {\Cal
O}_{{\Bbb  P}^1}(-2m+1)) \cong  {\Bbb C}  \oplus  {\Bbb C} \endaligned
\tag{$\ast\ast$} $$ Similarly, $X  +  L^n$  is  linearly   equivalent
to  $Q \cdot F$, where  $F$  is  a  hypersurface  of  degree  $a$.  As
in  the  case  II,  $X$  is projectively  Cohen-Macaulay. This
completes the  proof of Theorem  2.1.
\qed
\enddemo

\remark{Remark 2.6} For smooth surfaces in $\Bbb P^4$, Theorem 2.1 was
also proved by A.Aure and L.Roth , see \cite{Au} \cite{Ro}.
\endremark

\proclaim{Corollary 2.7}
Let $X^n  \subset \Bbb  P^{n+2}\,  (n \geq 2)$  be  a  locally
Cohen-Macaulay subvariety contained  in  a hyperquadric  $Q$ of rank
$3$(resp., $4$) with $\Sing{(Q)}=\Lambda$  Then,
\roster
\item[1]
$\depth_{\Lambda}{( {\Cal  I}_{X / Q} )}= 2$ (resp., $3$)
\item[2] If $X$ is not a complete intersection, then
$\Sing{(Q)}=\Lambda\subset X$.
\endroster
\endproclaim

\demo{Proof} For a proof of (1), we first assume that $\rank{(Q)}=3$.
Since $X$ is projectively Cohen-Macaulay, ${H}^i(\Bbb P^{n+2},{\Cal
I}_{X})={H}^i(Q,{\Cal  I}_{X / Q}) = 0$ for $i = 1,2.$ Therefore, by
(2.3.1), for $k=0$ or $k=1$, $$ H^1(\Cal U,{\Cal  I}_{X / Q}|_{\Cal
U}) \simeq \bigoplus_{m=0}^{\infty} S^m(V) \otimes  H^1(\Bbb P^1,
{\Cal  O}_{{\Bbb P}^1}(-2m-2a+k)) \simeq  H^2_{\Lambda}(Q,{\Cal  I}_{X
/ Q}), $$ which is obvious nonzero and so by Theorem 1.2.,
$\depth_{\Lambda}{( {\Cal  I}_{X / Q} )}= 2$. Similarly, in the case
of $\rank{(Q)}=4$, we have an exact sequence $$ 0 \rightarrow
{H}^2({\Cal U},{\Cal  I}_{X / Q}|_{\Cal  U}) \rightarrow
{H}^3_{\Lambda}(Q,{\Cal  I}_{X / Q})  \rightarrow {H}^3(Q,{\Cal  I}_{X
/ Q})  \rightarrow 0 $$ Note that for surface $X$, $h^3(Q, {\Cal I}_{X
/ Q})=h^3(\Bbb P^4, {\Cal I}_{X / {\Bbb P}^4}) =h^2(X, {\Cal O}_X)$
may not be zero. By the way, $\depth_{\Lambda}{( {\Cal  I}_{X / Q} )}=
3$ because $H^2(\Cal U,{\Cal I}_{X / Q} |_{\Cal  U}) \simeq
\bigoplus_{m=0}^{\infty}  S^m(V) \otimes  H^2(Q_1,  {\Cal  O}_{Q_1}
(-m-a,-m-b)),|a-b|  \leq 1 $ is nonzero.

For (2), by Theorem 2.1, we know that $X=Q \cap F_1 \cap F_2$, where
$F_i$ is a homogeneous polynomial of degree $a$, $\deg{(X)}=2a-1$ and
$Q \cap  F_i= X  \cup  L_i,$ \newline $L_i= < {\Lambda}, \ell_i >$
where $\ell_i$ is a line in $Q_1 \subset \Bbb P^3$ ( or $\ell_i$ is
one point of a smooth conic $C$ in $\Bbb P^2$ if $\rank{(Q)}=3$).
Thus, $\Lambda \subset X=Q \cap F_1 \cap F_2$ \qed
\enddemo

\remark{Remark 2.8}
It is well-known that if $X^n \subset {\Bbb P}^{n+2},n \geq 5$ and $X$
a locally Cohen-Macaulay variety, then $X$  is  a  complete
intersection if and only if $X$ is projectively Cohen-Macaulay, see
\cite {PSz}.
\endremark

\proclaim{Proposition 2.9}
If a locally Cohen-Macaulay variety  $X^n \subset {\Bbb P}^{n+2},(n
\ge 2)$ is contained in a hyperquadric $Q$ and is  not  a  complete
intersection, then \roster \item[1]  The saturated ideal $I_X  =
\bigoplus_{m  \in  {\Bbb  Z}}  H^0({\Bbb  P}^{n+2},  {\Cal  I}_X(m)) $
is generated by polynomials $Q,F_1,F_2,$  $\deg{(F_i)} = a$ and
$\deg{(X)}= 2a-1$.

\item[2]  
$\reg{(X)} = \frac{\deg{(X)}+1}{2}= a.$
\endroster
\endproclaim

\demo{Proof} Since $X$ is not a complete intersection, $\rank{(Q)}=3$
or $4$, as we have already noted. Assume that ${\Cal  I}_{X  /
Q}|_{\Cal  U}  =  {\pi}_{\Lambda}^{\ast}({\Cal   O}_{Q_1} (-a, -a+1))$
or ${\pi}_{\Lambda}^{\ast}(  i^{\ast}{\Cal  O}_{{\Bbb  P}^1}(-2a+1))$
as in the proof of Theorem 2.1. By ($\ast$) and ($\ast\ast$), we can
show that for $a \ge 2 $ $$ h^0(Q,  {\Cal I}_{X /Q}(k))  = \left\{
\aligned &0 \  \,\,\quad when \  2\le k < a,   \\ &2 \  \,\,\quad when
\  k = a. \endaligned \right. $$ Let $\alpha(I_X)= \min{\{i \in
\Bbb Z : h^0({\Cal I}_X(i)) \neq 0\}=2}$ and $\nu(I_X)=$number of
minimal generators of the saturated ideal $I_X$ of $X$. Since $X$ is
projectively Cohen-Macaulay of codimension two, we have  $\nu(I_X)=
\alpha(I_X) +1=3$ (see \cite{Mi}, Corollary 2.2.3.). Therefore, the
saturated ideal $I_X = \bigoplus_{m  \in  {\Bbb  Z}}  H^0({\Bbb
P}^{n+2},  {\Cal  I}_X(m)) $ is generated by polynomials $Q,F_1,F_2,$
$\deg{(F_i)} = a$. For a proof of (2), from ($\ast$) and ($\ast\ast$),
we have $h^0(Q, {\Cal I}_{X /Q}(a+1))  = 2n+4.$  The vector space of
all polynomials of the form $c_1F_1\cdot L_1 + c_2F_2\cdot L_2$, where
$L_1$ and $L_2$ are arbitrary linear forms, has dimension $(2n+6)$.
This implies that there are two relations among minimal generators
$Q,F_1,F_2$ in degree $(a+1).$ Thus, the minimal free resolution for
$I_X$ is a short exact sequence  as follows: $$ 0  \rightarrow S(-a-1)
\oplus S(-a-1)  \rightarrow S(-2) \oplus  S(-a) \oplus S(-a)
\rightarrow I_X  \rightarrow 0. $$ Consequently, $\reg{(X)} =
\frac{\deg{(X)}+1}{2}= a.$ \qed
\enddemo

\remark{Remark 2.10} If $X \subset \Bbb P^N $ is a complete
intersection defined by polynomials $F_1,F_2, \cdots , F_e$,
$\deg{(F_i)}=d_i$, $e=\codim{(X,\Bbb P^N)}$, then it is easily checked
that $\reg{(X)}=d_1+d_2+ \cdots +d_e-e+1.$
\endremark

\proclaim{Theorem  2.11}
Let  $X$  be  a  smooth  threefold in  ${\Bbb  P}^5$ which is  not
contained  in any  hyperquadric $Q$.  Then $\reg{(X)}   \leq
\deg{(X)}-3$. We also have $\reg{(S)} \leq \deg{(S)}-3$ for any
regular surface $S$ with $h^0({\Cal I}_{S / \Bbb P^4}(2))=0$ in
$\Bbb P^4$ except for the Veronese surface.
\endproclaim

\demo{Proof} 
The  following  construction  and  technique  are  basically due  to
Lazarsfeld  \cite{L} and Alzati  \cite{Al}. The author wants to give
more simplified and self-contained proof. Let ${\Bbb P}^5= \Proj{{\Bbb
C}[T_0,T_1,T_2,T_3,T_4,T_5]}$. Consider  the following diagram
$$
\CD
Bl_p({\Bbb  P}^5)  =  
{\Bbb  P}({\Cal  O}_{{\Bbb  P}^4}(1)  \oplus{\Cal  O}_{{\Bbb  P}^4})
@>  {\pi}_2  >>  
{\Bbb  P}^4  \\  
@V  {\pi}_1  VV  @.  @.  \\
X  \subset  {\Bbb  P}^5   @. 
\endCD
$$
where  $p = (0,0,0,0,0,1)  \notin  X$ and ${\pi}_p : {\Bbb  P}^5
\dashrightarrow  {\Bbb  P}^4$ is  a  generic  projection,  which
commutes with the  above  diagram,  i.e. ${\pi}_2  \circ  {\pi}^{-1}_1
\,  |_X  \,\,  =  \,\,  {\pi}_p  \,  |_X$. Therefore,  the  natural
morphism ${\Cal  O}_{{\Bbb  P}^5}(k)  \rightarrow  {\Cal  O}_X(k)$
induces  a  morphism  \newline $ w_k  : {\pi}_{2_{\ast}}(
{\pi}^{\ast}_1{\Cal  O}_{{\Bbb  P}^5}(k) )  \cong Sym^k( {\Cal
O}_{{\Bbb  P}^4}(1)  \oplus {\Cal  O}_{{\Bbb  P}^4})     \rightarrow
{\pi}_{2_{\ast}}({\pi}^{\ast}_1 {\Cal  O}_X(k) )  \simeq
{\pi}_{p_{\ast}}{\Cal  O}_X(k).$

By Theorem 1.4 (the  (dimension+2)-secant  Lemma), the  length  of
each  fibre  of ${\pi}_p$  is  at  most  4. Let $L_y$ be the line
through $p$,$y$. Then $$ w_3  \otimes  {\Bbb  C}(y)  : H^0(L_y,  {\Cal
O}_{L_y}(3))  \rightarrow H^0(L_ y,  {\Cal  O}_{{\pi}_p^{-1}(y)}(3))
$$ is  surjective, i.e., the scheme-theoretic fibre
${{\pi}_p^{-1}(y)}$ is 3-normal for  all  $y  \in  {\pi}_p(X)$ and by
tensoring  ${\Cal  O}_{{\Bbb  P}^4}(-3)$, we  have  an  exact
sequence, $$ 0  \rightarrow E  \rightarrow {\Cal  O}_{{\Bbb  P}^4}(-3)
\oplus {\Cal  O}_{{\Bbb  P}^4}(-2)  \oplus {\Cal  O}_{{\Bbb  P}^4}(-1)
\oplus {\Cal  O}_{{\Bbb  P}^4}  @>  w_3 \otimes {\Cal  O}_{{\Bbb
P}^4}(-3)  >> {\pi}_{p_{\ast}}{\Cal  O}_X  \rightarrow 0 $$ where  $E
=  \Ker{w_3}  \otimes {\Cal  O}_{{\Bbb  P}^4}(-3)$.

\proclaim{Lemma  2.12} The vector bundle $E = \Ker{w_3} \otimes
{\Cal  O}_{{\Bbb  P}^4}(-3)$  controls the  regularity  of  $X$,
i.e., $\reg{(X)}  \leq \reg{(E)}$.
\endproclaim

\demo{Proof}
According to Lazarsfeld \cite{L}, this is due to Gruson and Peskine.
Suppose  that  $E$  is  $(m+1)$-regular:  in other  words, $H^i({\Bbb
P}^4,  E(m+1-i))  =  0,  \,  i>0$. We'll  show  that  $X$  is
$(m+1)$-regular,  equivalently,  $X$  is  $m$-normal  and $H^i(X,
{\Cal  O}_X(m-i))  =  0,  \,\,  i \ge 1$. Consider the  sequence $$ 0
\rightarrow E  \rightarrow {\Cal  O}_{{\Bbb  P}^4}(-3)  \oplus {\Cal
O}_{{\Bbb  P}^4}(-2)  \oplus {\Cal  O}_{{\Bbb  P}^4}(-1)  \oplus {\Cal
O}_{{\Bbb  P}^4}  @>  w_3 \otimes {\Cal  O}_{{\Bbb  P}^4}(-3)  >>
{\pi}_{p_{\ast}}{\Cal  O}_X  \rightarrow 0 $$ Since  $H^1({\Bbb  P}^4,
E(m))  =  0$,   the  following  morphism $$ \CD H^0({\Cal  O}_{{\Bbb
P}^4}(m-3))  \oplus H^0({\Cal  O}_{{\Bbb  P}^4}(m-2))  \oplus
H^0({\Cal  O}_{{\Bbb  P}^4}(m-1))  \oplus H^0({\Cal  O}_{{\Bbb
P}^4}(m))  \\ @V  H^0(w_3)  VV  \\ H^0({\pi}_{p_{\ast}}{\Cal  O}_X(m))
\endCD $$ is  surjective,  where $ H^0({\pi}_{p_{\ast}}{\Cal  O}_X(m))
\simeq H^0(X, {\Cal  O}_X(m)) $ and $H^o(w_3)$  maps  \linebreak  $F_i
\in  H^0({\Cal  O}_{{\Bbb  P}^4}(m-i))$  to $F_i  \cdot  T_5^i$  for
$i  =  0,1,2,3$.  Hence  $X$  is  $m$-normal.

Secondly,  $H^i({\Bbb  P}^4,  E(m+1-i))  =  0,  \,\,  i \ge 2$
implies $H^i(X, {\Cal  O}_X(m-i))  =  0$  for  $i \ge 1$. Indeed,
from  the  exact  sequence $$ 0  \rightarrow E  \rightarrow {\Cal
O}_{{\Bbb  P}^4}(-3)  \oplus {\Cal  O}_{{\Bbb  P}^4}(-2)  \oplus {\Cal
O}_{{\Bbb  P}^4}(-1)  \oplus {\Cal  O}_{{\Bbb  P}^4}  @>  w_3  >>
{\pi}_{p_{\ast}}{\Cal  O}_X  \rightarrow 0 $$ we  get , for $ i \ge 1
$, $$ H^i(X, {\Cal  O}_X(m-i))  = H^i({\Bbb  P}^4,
{\pi}_{p_{\ast}}{\Cal  O}_X(m-i))  = H^{i+1}({\Bbb  P}^4,  E(m-i))  =
0 $$ Therefore,  $X$  is  $(m+1)$-regular, and  $\reg{(X)}  \leq
\reg{(E)}$. \qed
\enddemo
    
\proclaim{Claim}  $\text{E}^{\ast}$  is  $(-3)$-regular for a smooth
$3$-fold $X$ in ${\Bbb P}^5$.
\endproclaim

\demo{Proof  of Claim}
Note that  $E^{\ast}$ is  $(-3)$-regular if and only if $ H^i({\Bbb
P}^4,  \text{E}^{\ast}(-3-i))  =  0,  i = 1,2,3,4. $ For  $i = 1$, $
h^1(\text{E}^{\ast}(-4))  = h^3(\text{E}(-1))  =
h^2({\pi}_{p_{\ast}}{\Cal  O}_X(-1))  = h^2( {\Cal  O}_X(-1))  =  0 $
by  Kodaira's  vanishing  theorem.  When  $i  =  2$, $ h^2({\Bbb
P}^4,  \text {E}^{\ast}(-5))\allowmathbreak= h^2(\text E)=
h^1({\pi}_{p_{\ast}}{\Cal  O}_X)= h^1({\Cal  O}_X)=0 $ by  the
Barth-Larsen theorem \cite{B}.  For $i=3$, $ h^3({\Bbb  P}^4,  \text
{E}^{\ast}(-6))  = h^1(\text{E}(1))  =  0 $ because  we  have  the
following  exact  sequence $$ H^0( {\Cal  O}_{{\Bbb  P}^4}(-2)  \oplus
{\Cal  O}_{{\Bbb  P}^4}(-1)  \oplus {\Cal  O}_{{\Bbb  P}^4}  \oplus
{\Cal  O}_{{\Bbb  P}^4}(1))  @>  \alpha  >> H^0( {\Cal  O}_X(1))
\rightarrow H^1(\text {E}(1))  \rightarrow 0 $$ and  $\alpha$  is
surjective by Zak's linear normality theorem \cite{Z}. Finally,  if
$i=4$,  then $h^4(\text {E}^{\ast}(-7))  = h^0(\text {E}(2))  =  0$
because  $X$  is  not  contained  in any hyperquadric and  hence $$
H^0({\Bbb  P}^4, {\Cal  O}_{{\Bbb  P}^4}(-1)  \oplus {\Cal  O}_{{\Bbb
P}^4}  \oplus {\Cal  O}_{{\Bbb  P}^4}(1)  \oplus {\Cal  O}_{{\Bbb
P}^4}(2))  \rightarrow H^0({\Bbb  P}^4,  {\pi}_{p_{\ast}}{\Cal
O}_X(2)) $$ is  injective. Therefore, $E^{\ast}$  is $(-3)$-regular.
This  completes  the  proof  of  our  claim. \qed
\enddemo

Now,  let's  go  back  to  the  proof  of  Theorem  2.11. Note  that
$E  \simeq {\wedge}^3  E^{\ast}  \otimes  \det{E},
\rank{(E)} =  4 $. So,  by Lemma 1.1, $ \reg{(E)}  \leq
\reg{(\bigwedge^3 (E^{\ast}))}  + \reg{(\det{(E)})}\allowmathbreak
\leq d-3 $ where  $\operatorname{det}(E)  =  {\Cal  O}_{{\Bbb  P}^4}(-6-d)$
and $\reg{(\operatorname{det}(E))}  =  d+6$.  Thus $$ \reg{({\Cal
I}_X)}=\operatorname{reg}{(X)} \leq  \operatorname{reg}{(E)}  \leq
d-3 $$ Furthermore, suppose  that  $S$  is  a  smooth  regular
surface  in ${\Bbb  P}^4$  with  degree  $d$  and $h^0({\Cal  I}_S(2))
= 0$. By the same generic projection method, we  have  an  exact
sequence $$ 0  \rightarrow E  \rightarrow {\Cal  O}_{{\Bbb  P}^3}(-2)
\oplus {\Cal  O}_{{\Bbb  P}^3}(-1)  \oplus {\Cal  O}_{{\Bbb  P}^3}  @>
w_2 \otimes {\Cal  O}_{{\Bbb  P}^3}(-2) >> {\pi}_{p_{\ast}}{\Cal  O}_S
\rightarrow  0 $$ Since  $S$  is  linearly  normal and  $h^1({\Cal
O}_S)  =  0$, by  a  similar  method,  $E^{\ast}$  is $(-3)$-regular
and $ \operatorname{reg}{(S)} \leq \operatorname{reg}{(E)}  =
\operatorname{reg}{(\bigwedge^2 (E^{\ast})} \otimes
\operatorname{det}(E))  \leq  d-3$ \qed
\enddemo

\remark{\bf Remark  2.13}
\roster
\item"(a)" 
Such  a  bound  $(d-3)$  in  Theorem  2.11  is  the best  possible
because, as  shown  in example 2.0., the  dependency  locus  $X_0$  of
four  generic  sections  of ${\Omega}_{{\Bbb  P}^5}(2)$  has  degree
$7$  and  $\reg{(X)}  =4$. Such  a  threefold is  called  the
Palatini scroll  of  degree  $7$ which is not 3-regular. \item"(b)" By
Proposition 2.9 and Theorem 2.11, for a smooth 3-fold $X$ in $\Bbb
P^5$, $\reg{(X)}  \leq  d-1$ and $\reg{(X)}=d-1$ if and only if $X$
has  degree $3$ or $4$, i.e., it is either the Segre threefold
or the complete intersection of two quadrics.
\item"(c)"
Theorem 2.11 is also true for smooth 3-folds in $\Bbb P^5$ defined
over any algebraically closed field of characteristic zero because
Kodaira's vanishing theorem can be extended to this case. \endroster
\endremark

\proclaim{Corollary  2.14} Let  $X$  be  a  smooth  3-fold
in  ${\Bbb  P}^5$  of  degree  $d$.  Then $X$  is  k-normal  for  $k
\geq  d-4$.  The Palatini scroll of degree $7$ which, as we have
already noted,  is not 2-normal, shows that this is the optimal bound.
\endproclaim

\demo{Proof}
By Theorem  2.1  and  Theorem  2.11,  this  is  clear.   \qed
\enddemo

\bigskip

\enddocument